\documentclass[12pt]{Manuscript_Template_ICCA9_LaTex}

\usepackage[noadjust]{cite}
\usepackage{latexsym}
\usepackage{amsfonts}
\usepackage{amssymb}
\usepackage{amsthm}
\usepackage{amsmath}              
\usepackage{wasysym}              

\newcommand{\cl}{C \kern -0.1em \ell}     

\newcommand{\Mat}{{\rm Mat}}

  \newcommand{\be}{{\bf e}}

\newcommand{\bL}{{\bf L}}

  \newcommand{\bx}{{\bf x}}
  \newcommand{\by}{{\bf y}}

\newcommand{\Id}{{\bf 1}}

\newcommand{\BK}{\mathbb{K}}

\newcommand{\BC}{\mathbb{C}}
\newcommand{\BR}{\mathbb{R}}
\newcommand{\BH}{\mathbb{H}}
\newcommand{\BZ}{\mathbb{Z}}

\newcommand{\spn}{\mbox{\rm span}}

\newcommand{\beq}{\begin{equation}}
\newcommand{\eeq}{\end{equation}}

\input{diagxy}    

\newcommand{\ed}{\end{document}}
\def\ve{\varepsilon}
\def\vp{\varphi}
\def\Sg{{\hat{S}}}
\def\fg{\hat{f}}
\def\psig{{\psi_g}}
\def\phig{{\phi_g}}
\def\chpsi{\check{\psi}}
\def\chphi{\check{\phi}}

\newcommand{\ta}[2]{#1_#2\tilde{\phantom{.}}}
\newcommand{\cb}[1]{\mathcal{#1}}
\newcommand{\matrepr}[1]{\left[#1 \right]}
\newcommand{\Gpqe}[2]{G_{#1,#2}^\varepsilon}
\newcommand{\Gpq}[2]{G_{#1,#2}}
\newcommand{\Gpqf}[3]{G_{#1,#2}(#3)}

\newcommand{\Kpqf}[3]{K_{#1,#2}(#3)}
\newcommand{\Tpqf}[3]{T_{#1,#2}(#3)}

\newcommand{\iu}{{\underline{i}}}
\newcommand{\ju}{{\underline{j}}}

\newcommand{\fpower}[1]{{}^2 \kern -0.01em #1}

\newcommand{\chS}{\check{S}}
\newcommand{\chK}{\check{\BK}}
\renewcommand{\Id}{\mathrm{Id}}

\newcommand{\tve}{t_\ve}
\newcommand{\tp}{\ta{T}{\ve}}

\newcommand{\lang}{\mathopen{<}}
\newcommand{\rang}{\mathclose{>}}

\DeclareMathOperator{\hotimes}{\Hat{\otimes}}

\theoremstyle{plain}
\newtheorem{theorem}{Theorem}

\newtheorem{corollary}{Corollary}

\newtheorem{proposition}{Proposition}
\theoremstyle{definition}
\newtheorem{definition}{Definition}

\title{TRANSPOSITION ANTI-INVOLUTION IN CLIFFORD ALGEBRAS
AND INVARIANCE GROUPS OF SCALAR PRODUCTS ON SPINOR SPACES}

\author{\underline{R. \ Ab\l amowicz}$^*$  and B. Fauser$^+$}

\address{\parbox[t]{0.45\textwidth}{%
    \hskip-1.1ex$^*$Department of Mathematics,\\
    Tennessee Technological University\\
    Cookeville, TN 38505, USA\\
    \selectfont \normalfont E-mail: rablamowicz@tntech.edu}
\parbox[t]{0.45\textwidth}{%
    \hskip-1.15ex$^+$School of Computer Science,\\
    The University of Birmingham\\
    Edgbaston-Birmingham, B15 2TT, UK\\
    \selectfont \normalfont E-mail: b.fauser@cs.bham.ac.uk}
}


\abstract{We introduce on the abstract level in real Clifford algebras
$\cl_{p,q}$ of a non-degenerate quadratic space $(V,Q)$, where $Q$ has
signature $\varepsilon=(p,q)$, a transposition anti-involution $\tp$. In
a spinor representation, the anti-involution $\tp$ gives transposition,
complex Hermitian conjugation or quaternionic Hermitian conjugation when
the spinor space $\check{S}$ is viewed as a $\cl_{p,q}$-left and
$\check{\BK}$-right module with $\check{\BK}$ isomorphic to $\BR$ or
$\fpower\BR$, $\BC$, or, $\BH$ or $\fpower \BH$. This map and its
application to SVD was first presented at ICCA~7 in Toulouse in
2005~\cite{ablamowicz2005}. 

The anti-involution $\tp$ is a lifting to $\cl_{p,q}$ of an orthogonal
involution $\tve: V \rightarrow V$ which depends on the signature of $Q$.
The involution is a symmetric correlation~\cite{porteous}
$\tve: V \rightarrow V^{*} \cong V$ and it allows one to define a
reciprocal basis for the dual space $(V^{*},Q)$. When the Clifford algebra
$\cl_{p,q}$ splits into the graded tensor product
$\cl_{p,0} \hotimes \cl_{0,q}$, the anti-involution $\tp$ acts as reversion
on $\cl_{p,0}$ and as conjugation on $\cl_{0,q}$. Using the concept of a
transpose of a linear mapping one can show that if $[L_u]$ is a matrix in
the left regular representation of the operator
$L_u: \cl_{p,q} \rightarrow \cl_{p,q}$ relative to a Grassmann basis
$\cb{B}$ in $\cl_{p,q},$ then matrix $[L_{\tp(u)}]$ is the matrix transpose
of $[L_u]$, see~\cite{part1}. 

Of particular importance is the action of $\tp$ on the spinor space. The
algebraic spinor space~$\check{S}$ is realized as a left minimal ideal
generated by a primitive idempotent $f$, or a sum $f+\hat{f}$ in simple or
semisimple algebras as in~\cite{lounesto}. The map $\tp$ allows us to define
a new spinor scalar product
$\check{S} \times \check{S} \rightarrow \check{\BK}$, where
$\BK=f\cl_{p,q}f$ and $\check{\BK}=\BK$ or $\BK \oplus \hat{\BK}$ depending
whether the algebra is simple or semisimple. Our scalar product is in general
different from the two scalar products discussed in literature,
e.g.,~\cite{lounesto}. However, it reduces to one or the other in Euclidean
and anti-Euclidean signatures. The anti-involution $\tp$ acts as the identity
map, complex conjugation, or quaternionic conjugation on $\check{\BK}.$ Thus,
the action of $\tp$ on spinors results in matrix transposition, complex
Hermitian conjugation, or quaternionic Hermitian conjugation. We classify
automorphism group of the new product as
$O(N)$, $U(N)$, $Sp(N)$, $\fpower O(N)$, or $\fpower Sp(N)$.
}

\begin{document}
\section{Introduction}
Let $\cl_n$ be a universal Clifford algebra over an $n$-dimensional real quadratic space $(V,Q)$ with 
$Q(\bx) = \ve_1x_1^2 + \ve_2x_2^2 + \cdots + \ve_nx_n^2$ where $\ve_i = \pm 1$ and 
$\bx = x_1\be_1 + \cdots + x_n\be_n \in V$ for an orthonormal basis
$\mathcal{B}_1 = \{\be_i\}_{i=1}^{n}$. Let $\cb{B}$ be the canonical basis of 
$\bigwedge V$ generated by $\cb{B}_1$. That is, let $[n]=\{1,2,\ldots,n\}$ and
denote arbitrary, canonically ordered subsets of $[n]$, by underlined Roman
characters. The basis elements of $\bigwedge V$, or, of $\cl_n$ due to the
linear space isomorphism $\bigwedge V \rightarrow \cl_n$~\cite{lounesto}, can
be indexed by these finite ordered subsets as
$\be_\iu = \wedge_{i \in \iu}\, \be_i$. Then, an arbitrary element of
$\bigwedge V \cong \cl_n$ can be written as
$u = \sum_{\iu \in 2^{[n]}} u_\iu \be_{\iu}$ where $u_\iu \in \BR$ for each
$\iu \in 2^{[n]}$. The unit element $1$ of $\cl_n$ is identified with
$\be_\emptyset$. Our preferred basis for $\cl_n$ is the exterior algebra
basis $\cb{B}$ sorted by an \textit{admissible} monomial order~$\prec$ on
$\bigwedge V.$ We choose for $\prec$ the monomial order called
$\mathtt{InvLex}$, or, the \textit{inverse lexicographic order}~\cite{ablamowicz2009a,GfG}. Let $B$ be
the symmetric bilinear form defined by $Q$ and let
$\lang \cdot,\cdot \rang: \bigwedge V \times \bigwedge V \rightarrow \BR$ be
an extension of $B$ to $\bigwedge V$~\cite{lounesto}. We will need this
extension later when we define the Clifford algebra $\cl(V^\ast,Q)$.

We begin by defining the following map on $(V,Q)$ dependent on the signature
$\ve$ of~$Q$.
\begin{definition}
Let $t_\ve: V \rightarrow V$ be the linear map defined as
\begin{equation}
   t_\ve(\bx) = t_\ve(\sum_{i=1}^n x_i \be_i)
              = \sum_{i=1}^n x_i \left(\frac{\be_i}{\ve_i}\right)
	      = \sum_{i=1}^n x_i \left(\ve_i \be_i\right)
\label{eq:tve} 
\end{equation}
for any $\bx \in V$ and for the orthonormal basis
$\cb{B}_1 = \{\be_i\}_{i=1}^{n}$ in $V$ diagonalizing $Q$.
\end{definition}
The $t_\ve$ map can be viewed in two ways: (1) As a linear orthogonal
involution of $V$; (2) As a \textit{correlation}~\cite{porteous} mapping
$t_\ve: V \rightarrow V^\ast \cong V$. The set of vectors 
$\cb{B}_1^\ast = \{t_\ve(\be_i)\}_{i=1}^{n}$ gives an orthonormal basis
in the dual space $(V^\ast,Q)$. Furthermore, under the identification
$V \cong V^\ast$, $t_\ve$ is a symmetric non-degenerate correlation
on~$V$ thus making the pair $(V,t_\ve)$ into a
\textit{non-degenerate real correlated (linear) space}~\cite{part1}. Then,
viewing $\tve$ as a correlation $V \rightarrow V^\ast$, we can define the
action of $t_\ve(\bx)\in V^\ast$ on $\by \in V$ for any $\bx \in V$ as
\begin{equation}
t_\ve(\bx)(\by) = \lang t_\ve(\bx), \by \rang,
\label{eq:action1}
\end{equation}
and we get the expected duality relation among the basis elements in
$\cb{B}_1$ and $\cb{B}_1^\ast$:
\begin{equation}
t_\ve(\be_i)(\be_j) = \lang \ve_i \be_i, \be_j  \rang = 
\ve_i \lang \be_i, \be_j \rang = \ve_i\ve_j\delta_{i,j} = \delta_{i,j}.
\label{eq:action2}
\end{equation}
The extension of the duality $V \rightarrow V^\ast$ to the Clifford algebras
$\cl(V,Q) \rightarrow \cl(V^\ast,Q)$ is of fundamental importance to defining
a new transposition scalar product on spinor spaces. When we apply Porteous'
theorem~\cite[Thm. 15.32]{porteous} to the involution $\tve$, we get the
following theorem and its corollary proven
in~\cite{part1}.\footnote{%
We view $\cl_n$ as Porteous' $\bL^\alpha$-\textit{Clifford algebra for $(V,Q)$}
under the identification $\bL=\BR$ and $\alpha = 1_\BR$.}
\begin{proposition}
Let $\cb{A}=\cl_n$ be the universal Clifford algebra of $(V,Q)$ and let
$t_\ve:V \rightarrow V$ be the orthogonal involution of $V$ defined
in~(\ref{eq:tve}). Then there exists a unique algebra involution
$T_\ve$ of $\cb{A}$ and a unique algebra anti-involution $\ta{T}{\ve}$ of
$\cb{A}$ such that the following diagrams commute:
\begin{equation}
\bfig
\square|arra|/>` >->` >->`>/[V`V`\cb{A}`\cb{A};t_\ve`\iota`\iota`T_\ve]
\efig
\qquad \mbox{and} \qquad
\bfig
\square|arra|/>` >->` >->`>/[V`V`\cb{A}`\cb{A};t_\ve`\iota`\iota`\ta{T}{\ve}]
\efig
\label{d:diag1}
\end{equation}
In particular, we can define $T_\ve$ and $\ta{T}{\ve}$ as follows:
\begin{itemize}
\item[(i)] For simple $k$-vectors $\be_\iu$ in $\cb{B}$, let 
     $T_\ve(\be_\iu) = T_\ve (\prod_{i \in \iu} \be_i) 
                     = \prod_{i \in \iu} t_\ve(\be_i)$ where
     $k= |\iu|$ and $T_\ve(1_\cb{A}) = 1_\cb{A}$. Then, extend by linearity
     to all of $\cb{A}$.
\item[(ii)] For simple $k$-vectors $\be_\iu$ in $\cb{B}$, let 
  \begin{gather}
     \ta{T}{\ve}(\be_\iu) 
        = \ta{T}{\ve} (\prod_{i \in \iu} \be_i)
        = (\prod_{i \in \iu} t_\ve(\be_i))\tilde{}
        = (-1)^{\frac{k(k-1)}{2}}\prod_{i \in \iu} t_\ve(\be_i)
\label{eq:defT}
\end{gather} 
where $k=|\iu|$ and $\ta{T}{\ve}(1_\cb{A}) = 1_\cb{A}$. Then, extend by
linearity to all of $\cb{A}$.
\end{itemize}
\label{p:p1}
\end{proposition}
Maple code of the procedure~$\mathtt{tp}$ which implements the anti-involution
$\ta{T}{\ve}$ in~$\cl_n$, was first presented at ICCA 7 in
Toulouse~\cite{ablamowicz2005}. The procedure $\mathtt{tp}$ requires
the $\mathtt{CLIFFORD}$ package~\cite{clifford3}. In the following
corollary, $\alpha,\beta,\gamma$ denote, respectively, the grade involution,
the reversion, and the conjugation in~$\cl_n$.
\begin{corollary}
Let $\cb{A} = \cl_{p,q}$ and let $T_\ve: \cb{A} \rightarrow \cb{A}$ and
$\ta{T}{\ve}: \cb{A} \rightarrow \cb{A}$ be the involution and the
anti-involution of $\cb{A}$ from Proposition~\ref{p:p1}.
\begin{itemize}
\item[(i)] For the Euclidean signature $(p,q)=(n,0)$, or $p-q=n$, we have
   $t_\ve = 1_V$. Thus, $T_\ve$ is the identity map $1_\cb{A}$ on $\cb{A}$
   and $\ta{T}{\ve}$ is the reversion~$\beta$ of $\cb{A}$.
\item[(ii)] For the anti-Euclidean signature $(p,q)=(0,n)$, or $p-q=-n$, we
   have $t_\ve = -1_V$. Thus, $T_\ve$ is the grade involution $\alpha$ of
   $\cb{A}$ and $\ta{T}{\ve}$ is the conjugation~$\gamma$ of $\cb{A}$.
\item[(iii)] For all other signatures $-n < p-q < n$, we have
   $t_\ve = 1_{V_1} \otimes -1_{V_2}$ where
   $(V,Q) = (V_1,Q_1) \perp (V_2,Q_2)$. Here, $(V_1,Q_1)$ is the Euclidean
   subspace of $(V,Q)$ of dimension~$p$ spanned by $\{\be_i\}_{i=1}^{p}$ with
   $Q_1 =Q|_{V_1}$ while $(V_2,Q_2)$ is the anti-Euclidean subspace of $(V,Q)$
   of dimension~$q$ spanned by $\{\be_i\}_{i=p+1}^{n=p+q}$ with $Q_2 =Q|_{V_2}$.
   Let $\cb{A}_1 = \cl(V_1,Q_1)$ and $\cb{A}_2 = \cl(V_2,Q_2)$ so
   $\cl(V,Q) \cong \cl(V_1,Q_1) \hotimes \cl(V_2,Q_2)$. Let $S$ (resp.
   $\hat{S}$) be the \emph{ungraded switch} (resp. the \emph{graded switch})
   on $\cl(V_1,Q_1) \hotimes \cl(V_2,Q_2)$.\footnote{%
The switches are defined on the basis tensors 
$\be_\iu \hotimes \be_\ju \in \cl_{p,0} \hotimes \cl_{0,q}$ as 
$S(\be_\iu \hotimes \be_\ju) = \be_\ju \hotimes \be_\iu$ and 
$\hat{S}(\be_\iu \hotimes \be_\ju) = (-1)^{|\iu||\ju|}\be_\ju \hotimes \be_\iu$.
Then, their action is extended by linearity to the graded product 
$\cl_{p,0} \hotimes \cl_{0,q}$~\cite{part1}} 
   Then, 
   $$
   T_\ve = 1_{\cb{A}_1} \otimes \alpha_{\cb{A}_2} \quad \mbox{and} 
   \quad \ta{T}{\ve} = (\beta_{\cb{A}_1} \otimes\,  \gamma_{\cb{A}_2}) 
         \, \circ \,(\hat{S}\, \circ\, S).
   $$ 
\item[(iv)] The anti-involution $\ta{T}{\ve}$ is related to the involution
   $T_\ve$ through the reversion~$\beta$ as follows:
   $\ta{T}{\ve} = T_\ve \circ \beta =  \beta \circ T_\ve$.
\end{itemize}
\end{corollary}
For an extensive discussion of the properties of the involutions $\ta{T}{\ve}$
and $T_\ve$ see~\cite{part1}.

Since $(V^\ast,Q)$ is a non-degenerate quadratic space spanned by the
orthonormal basis $\cb{B}^\ast_1$, we can define the Clifford algebra
$\cl(V^\ast,Q)$ as expected. 
\begin{definition}
  The \textit{Clifford algebra over the dual space} $V^\ast$ is the universal
  Clifford algebra $\cl(V^\ast,Q)$ of the quadratic pair $(V^\ast,Q)$. For
  short, we denote this algebra by $\cl_n^\ast$.\footnote{%
Although from now on we denote the Clifford algebra of the dual $(V^\ast,Q)$
via $\cl_n^\ast$, we do not claim that $\cl_n^\ast$ is the \textit{dual}
algebra of $\cl_n$ in categorical sense as it was considered
in~\cite{oziewicz1997} and references therein.}
\label{def:dualcliff}
\end{definition}
Let $\cb{B}^\ast$ be the canonical basis of $\bigwedge V^\ast \cong \cl^\ast_n$
generated by $\cb{B}^\ast_1$ and sorted by~$\mathtt{InvLex}$. That is, we
define $\cb{B}^\ast = \{T_\ve(\be_\iu) \,| \, \be_\iu \in \cb{B}\}$ given that
\begin{equation}
  \lang T_\ve(\be_\iu),\be_\ju \rang = \delta_{\iu,\ju} 
  \label{eq:dualbasisClstar}
\end{equation}
for $\be_\iu,\be_\ju \in \cb{B}$ and $T_\ve(\be_\iu) \in \cb{B}^\ast$. An
arbitrary linear form $\vp$ in $\bigwedge V^\ast \cong \cl^\ast_n$ can be
written as
\begin{equation}
  \vp = \sum_{\iu \in 2^{[n]}} \vp_\iu T_\ve(\be_\iu)
  \label{eq:vp}
\end{equation}
where $\vp_\iu \in \BR$ for each $\iu \in 2^{[n]}$. Due to the linear
isomorphisms $V \cong V^\ast$ and $\bigwedge V^\ast \cong \cl(V^\ast,Q)$, we
extend, by a small abuse of notation, the inner product
$\lang \cdot,\cdot \rang$ defined in $\bigwedge V$ to 
\begin{equation}
  \lang \cdot,\cdot \rang: \bigwedge V^\ast \times \bigwedge V^\ast \rightarrow \BR.
  \label{eq:dualprod}
\end{equation}
In this way we find, as expected, that the matrix of this inner product on
$\bigwedge V^\ast$ is also diagonal, that is, that the basis~$\cb{B}^\ast$ is
orthonormal with respect to $\lang \cdot, \cdot \rang$. We extend the action
of dual vectors from $V^\ast$ on $V$ to all linear forms $\varphi$ in
$\cl^\ast_n$ acting on multivectors $v$ in $\cl_n$ via the inner
product~(\ref{eq:dualprod}) as
\begin{equation}
  \vp(v) = \lang \vp,v \rang
         = \sum_{\iu \in 2^{[n]}} \vp_\iu v_\iu
  \label{eq:phiaction}
\end{equation}
given that $\vp = \sum_{\iu \in 2^{[n]}} \vp_\iu T_\ve(\be_\iu) \in \cl^\ast_n$
where $\vp_\iu = \vp(\be_\iu) \in \BR$ and
$v = \sum_{\iu \in 2^{[n]}} v_\iu \be_\iu  \in \cl_n$ for some coefficients
$v_\iu \in \BR$.

Properties of the left multiplication operator $L_u: \cl_n \rightarrow \cl_n$,
$v \mapsto uv,\, \forall v \in \cl_n$ and its dual $L_{\tilde{u}}$ with respect
to the inner product
$\lang \cdot,\cdot \rang: \bigwedge V \times \bigwedge V \rightarrow \BR$ are
discussed in~\cite{part1}. In particular, it is shown there that if
$\matrepr{L_u}$ is the matrix of the operator $L_u$ relative to the basis
$\cb{B}$ and $\matrepr{L_{\ta{T}{\ve}(u)}}$ is the matrix of the operator
$L_{\ta{T}{\ve}(u)}$ relative to the basis $\cb{B}$, then 
$\matrepr{L_u}^T = \matrepr{L_{\ta{T}{\ve}(u)}} 
  = \matrepr{  L_{T_\ve(\tilde{u})}}$ where $\matrepr{L_u}^T$ is the matrix
transpose of $\matrepr{L_u}$. However, in order to introduce a new scalar
product on spinor spaces related to the involution $\tp$, we need to discuss
the action of $\tp$ on spinor spaces.

\section{Action of the transposition involution on spinor spaces}

Stabilizer groups $\Gpqf{p}{q}{f}$ of primitive idempotents $f$ are classified
in~\cite{part2}. The stabilizer $\Gpqf{p}{q}{f}$ is a normal subgroup of
Salingaros' finite vee group
$\Gpq{p}{q}$~\cite{salingaros1,salingaros2,salingaros3} which acts via
conjugation on $\cl_{p,q}$. The importance of the stabilizers to the spinor
representation theory lies in the fact that a \textit{transversal}\footnote{%
Let $K$ be a subgroup of a group $G$. A \textit{transversal} $\ell$ of $K$
in~$G$ is a subset of~$G$ consisting of exactly one element $\ell(bK)$ from
every (left) coset $bK$, and with $\ell(K)=1$~\cite{rotman}.}
of $\Gpqf{p}{q}{f}$ in $\Gpq{p}{q}$ generates spinor bases in $S=\cl_{p,q}f$
and $\Sg=\cl_{p,q}\fg$. In~\cite{part2} it is also shown that depending on the
signature $\varepsilon=(p,q)$, the real anti-involution $\tp$ is responsible
for transposition, the Hermitian complex, or the Hermitian quaternionic
conjugation of a matrix $[u]$ for any $u$ in all Clifford algebras $\cl_{p,q}$
with the spinor representation realized either in $S$ (simple algebras) or in
$\check{S}=S \oplus \Sg$ (semisimple algebras).  This is because $\tp$ acts
on $\BK = f\cl_{p,q}f$ and $\check{\BK}=\BK \oplus \hat{\BK}$ as an
anti-involution. Thus, $\tp$ allows us to define a dual spinor space $S^{*}$
or $\check{S}^{*}$, a new spinor product, and a new spinor norm. The following
results are proven in~\cite{part2}.
\begin{proposition}
Let $\psi,\phi \in S = \cl_{p,q}f$. Then, $\ta{T}{\ve}(\psi)\phi \in \BK$. In
particular, $\ta{T}{\ve}(\psi)\psi \in \BR f \subset \BK$.
\label{prop5}
\end{proposition}
Thus, we can define an invariance group of the scalar product
$S \times S \rightarrow \BK$, $(\psi,\phi) \mapsto \ta{T}{\ve}(\psi)\phi$, as
follows:
\begin{definition}
Let $G_{p,q}^\ve = \{ g \in \cl_{p,q}\; |\; \ta{T}{\ve}(g)g = 1\}$.
\end{definition}
\noindent
We find that $G_{p,q}(f) \unlhd G_{p,q} \leq G_{p,q}^\ve < \cl_{p,q}^{\times}$ (the group of units in $\cl_{p,q}$). Let $\cb{F}=\{f_i\}_{i=1}^N$ be a set of $N=2^k$, $k=q-r_{q-p}$, mutually annihilating primitive idempotents adding up to $1$ in a simple Clifford algebra 
$\cl_{p,q}$.\footnote{Here, $r_i$ is Radon-Hurwitz number defined by recursion as $r_{i+8}=r_i+4$ and these initial values: $r_0=0, r_1=1, r_2=r_3=2, r_4=r_5=r_6=r_7=3$~\cite{hahn,lounesto}.} The set $\cb{F}$ constitutes one orbit under the action of $G_{p,q}$~\cite{part2}.
\begin{proposition}
Let $\cl_{p,q}$ be a simple Clifford algebra, $p-q \neq 1 \bmod 4$ and
$p+q\leq 9$. Let $\psi_i \in S_i = \cl_{p,q}f_i$, $f_i \in \cb{F},$ and let $[\psi_i]$
(resp. $[\tp(\psi_i)]$) be the matrix of $\psi_i$ (resp. $\tp(\psi_i)$) in
the spinor representation with respect to the ordered basis
$\cb{S}_1 = [m_1 f_1,\ldots,m_N f_1]$ with $\alpha_i=m_i^2$.\footnote{%
For the sake of consistency with a proof of this proposition given
in~\cite{part2} we remark that $\alpha_i$ is just the square of the monomial
$m_i^2 \in \{\pm 1\}$.} Then, 
\begin{equation}
[\tp(\psi_i)] = \begin{cases} 
    [\psi_i]^T        & \textit{if $p-q =0,1,2 \bmod 8;$} \\
    [\psi_i]^\dagger  & \textit{if $p-q =3,7 \bmod 8;$} \\
    [\psi_i]^\ddagger & \textit{if $p-q =4,5,6 \bmod 8;$}
\end{cases}
\label{eq:daggers}
\end{equation}
where $T$ denotes transposition, $\dagger$ denotes Hermitian complex
conjugation, and $\ddagger$ denotes Hermitian quaternionic conjugation.
\label{prop6}
\end{proposition}
This action of $\tp$ on $S=S_i$ extends to a similar action on $\Sg$, hence to
$\check{S} = S \oplus \Sg$ as it is shown in~\cite{part2,part3}. In particular,
the product $(\psi,\phi) \mapsto \ta{T}{\ve}(\psi)\phi$ is invariant under
two of the subgroups of $G_{p,q}^\ve$: The Salingaros' vee group
$\Gpq{p}{q} < \Gpqe{p}{q}$ and the stabilizer group $\Gpqf{p}{q}{f}$ of a
primitive idempotent $f$. Since the stabilizer group $\Gpqf{p}{q}{f}$ and its
subgroups play an important role in constructing and understanding spinor
representation of Clifford algebras, we provide here a brief summary of
related definitions and findings. See~\cite{part3} for a complete
discussion.

Primitive idempotents $f \in \cb{F} \subset \cl_{p,q}$ formed out of commuting basis
monomials $\be_{\iu_1},\ldots,\be_{\iu_k}$ in~$\cb{B}$ with square~$1$ have
the form $f = \frac12(1\pm \be_{\iu_1})\frac12(1\pm \be_{\iu_2})\cdots 
\frac12(1\pm \be_{\iu_k})$ where $k = q - r_{q-p}$. With any primitive idempotent~$f$, we associate the following groups:\\
(i) The \textit{stabilizer $\Gpqf{p}{q}{f}$ of $f$} defined as
\begin{gather}
  \Gpqf{p}{q}{f} = \{ m \in \Gpq{p}{q} \mid m f m^{-1} = f \} < \Gpq{p}{q}.
\end{gather}
The stabilizer $\Gpqf{p}{q}{f}$ is a normal subgroup of $\Gpq{p}{q}$. In
particular,
\begin{gather}
 |\Gpqf{p}{q}{f}| = \begin{cases} 2^{1+p+r_{q-p}}, & p - q \neq 1 \bmod 4;\\
                                  2^{2+p+r_{q-p}}, & p - q = 1 \bmod 4.
                    \end{cases}
\label{eq:orderGpqf}
\end{gather}
(ii) An abelian \emph{idempotent group} $\Tpqf{p}{q}{f}$ \textit{of $f$}, a
subgroup of $\Gpqf{p}{q}{f}$ defined as
\begin{gather}
  \Tpqf{p}{q}{f} 
   = \langle \pm 1, \be_{\iu_1},\ldots, \be_{\iu_k} \rangle < \Gpqf{p}{q}{f},
\label{def:Tpqf}
\end{gather}
where $k= q - r_{q-p}.$\\

(iii) A \emph{field group} $\Kpqf{p}{q}{f}$ \textit{of $f$}, a subgroup of
$\Gpqf{p}{q}{f}$, related to the (skew double) field $\BK \cong f\cl_{p,q}f$,
and defined as
\begin{gather}
   \Kpqf{p}{q}{f} =  \langle \pm 1, m \mid m \in \cb{K}\rangle < \Gpqf{p}{q}{f}
   \label{def:Kpqf}
\end{gather}
where $\cb{K}$ is a set of Grassmann monomials in~$\cb{B}$ which provide a
basis for~$\BK=f\cl_{p,q}f$ as a real subalgebra of $\cl_{p,q}$.

The following theorem proven in~\cite{part3} relates the above groups to
$\Gpq{p}{q}$ and its commutator subgroup $\Gpq{p}{q}'$.\footnote{%
We have $\Gpq{p}{q}'=\{1,-1\}$ since any two monomials in $\Gpq{p}{q}$ either
commute or anticommute.}
\begin{theorem}
Let $f$ be a primitive idempotent in a simple or semisimple Clifford algebra
$\cl_{p,q}$ and let $\Gpq{p}{q}$, $\Gpqf{p}{q}{f}$, $\Tpqf{p}{q}{f}$,
$\Kpqf{p}{q}{f}$, and $\Gpq{p}{q}'$ be the groups defined above. Furthermore,
let $S=\cl_{p,q}f$ and $\BK=f\cl_{p,q}f$. 
\begin{itemize}
\setlength{\itemsep}{0.25ex}
\item[(i)] Elements of $\Tpqf{p}{q}{f}$ and $\Kpqf{p}{q}{f}$ commute.
\item[(ii)] $\Tpqf{p}{q}{f} \cap \Kpqf{p}{q}{f} = \Gpq{p}{q}' = \{\pm 1 \}$.
\item[(iii)] $\Gpqf{p}{q}{f} = \Tpqf{p}{q}{f}\Kpqf{p}{q}{f} = 
\Kpqf{p}{q}{f}\Tpqf{p}{q}{f}$.
\item[(iv)]  $|\Gpqf{p}{q}{f}| = |\Tpqf{p}{q}{f}\Kpqf{p}{q}{f}| =
\frac12 |\Tpqf{p}{q}{f}||\Kpqf{p}{q}{f}|$. 
\item[(v)] $\Gpqf{p}{q}{f} \lhd \Gpq{p}{q}$, $\Tpqf{p}{q}{f} \lhd \Gpq{p}{q}$,
and $\Kpqf{p}{q}{f} \lhd \Gpq{p}{q}$. In particular, $\Tpqf{p}{q}{f}$ and
$\Kpqf{p}{q}{f}$ are normal subgroups of $\Gpqf{p}{q}{f}$.
\item[(vi)] $\Gpqf{p}{q}{f} /\Kpqf{p}{q}{f} \cong \Tpqf{p}{q}{f} /\Gpq{p}{q}'$ and 
$\Gpqf{p}{q}{f} /\Tpqf{p}{q}{f} \cong \Kpqf{p}{q}{f} /\Gpq{p}{q}'$.
\item[(vii)] $(\Gpqf{p}{q}{f}/\Gpq{p}{q}')/(\Tpqf{p}{q}{f}/\Gpq{p}{q}')\cong \Gpqf{p}{q}{f}/\Tpqf{p}{q}{f} \cong \Kpqf{p}{q}{f}/\{\pm 1 \} $ and the transversal of $\Tpqf{p}{q}{f}$ in $\Gpqf{p}{q}{f}$ spans $\BK$ over $\BR$ modulo~$f$.
\item[(viii)] A transversal of $\Gpqf{p}{q}{f}$ in $\Gpq{p}{q}$ spans $S$ over $\BK$ modulo~$f$.  
\item[(ix)] $(\Gpqf{p}{q}{f}/\Tpqf{p}{q}{f}) \lhd
(\Gpq{p}{q}/\Tpqf{p}{q}{f})$ and 
$(\Gpq{p}{q}/\Tpqf{p}{q}{f})/(\Gpqf{p}{q}{f}/\Tpqf{p}{q}{f}) \cong \Gpq{p}{q}/\Gpqf{p}{q}{f}$
and a transversal of $\Tpqf{p}{q}{f}$ in $\Gpq{p}{q}$ spans $S$  over 
$\BR$ modulo~$f$.
\item[(x)] The stabilizer 
$\Gpqf{p}{q}{f} = \bigcap_{x \in \Tpqf{p}{q}{f}} C_{\Gpq{p}{q}}(x) = C_{\Gpq{p}{q}}(\Tpqf{p}{q}{f})$ where $C_{\Gpq{p}{q}}(x)$ is the centralizer of $x$ in $\Gpq{p}{q}$ and
$C_{\Gpq{p}{q}}(\Tpqf{p}{q}{f})$ is the centralizer of $\Tpqf{p}{q}{f}$ in
$\Gpq{p}{q}$.
\end{itemize}
\label{maintheorem}
\end{theorem}
Recall that in $\mathtt{CLIFFORD}$~\cite{clifford3} information about each
Clifford algebra $\cl_{p,q}$ for $p+q\le 9$ is stored in a built-in data file.
This information can be retrieved in the form of a seven-element list with the
command $\texttt{clidata([p,q])}$. For example, for $\cl_{3,0}$ we find:
\begin{gather}
  \mathtt{data} = [complex,
                   2,
		   simple,
		   \frac12 \Id+ \frac12 \be_1,
		   [\Id, \be_2, \be_3, \be_{23}],
		   [\Id, \be_{23}], [\Id, \be_{2}]]
  \label{eq:sampledata}
\end{gather}
where $\Id$ denotes the identity element of the algebra. In particular, from
the above we find that: (i) $\cl_{3,0}$ is a simple algebra isomorphic to
$\Mat(2,\BC)$; ($\mathtt{data[1]}$, $\mathtt{data[2]}$, $\mathtt{data[3]}$)
(ii) The expression $\frac12 \Id+ \frac12 \be_1$ ($\mathtt{data[4]}$) is a
primitive idempotent $f$ which may be used to generate a spinor ideal
$S=\cl_{3,0}f$;
(iii) The fifth entry $\mathtt{data[5]}$ provides, modulo~$f$, a real basis
for $S$, that is, $S = \spn_\BR \{ f, \be_2 f, \be_3 f, \be_{23} f \}$;
(iv) The sixth entry $\mathtt{data[6]}$ provides, modulo~$f$, a real basis for
$\BK=f \cl_{3,0}f \cong \BC$, that is, $\BK = \spn_\BR \{f, \be_{23}f\}$;
and, (v) The seventh entry $\mathtt{data[7]}$ provides, modulo~$f$, a basis
for $S$ over $\BK$, that is, $S = \spn_\BK \{ f, \be_2 f\}$.\footnote{%
See \cite{ablamowicz1996,ablamowicz1998} how to use \texttt{CLIFFORD}.}

The above theorem yields the following corollary:
\begin{corollary}
Let $\mathtt{data}$ be the list of data returned by the procedure
$\mathtt{clidata}$ in $\mathtt{CLIFFORD}$. Then, $\mathtt{data[5]}$ is a
transversal of $\Tpqf{p}{q}{f}$ in $\Gpq{p}{q}$; $\mathtt{data[6]}$ is a
transversal of $\Tpqf{p}{q}{f}$ in $\Gpqf{p}{q}{f}$; and $\mathtt{data[7]}$
is a transversal of $\Gpqf{p}{q}{f}$ in $\Gpq{p}{q}$. Therefore,
$|\mathtt{data[5]}| = |\mathtt{data[6]}||\mathtt{data[7]}|$. This is
equivalent to $\vert\frac{\Gpq{p}{q}}{\Tpqf{p}{q}{f}}\vert =
\vert\frac{\Gpqf{p}{q}{f}}{\Tpqf{p}{q}{f}}\vert\, 
\vert\frac{\Gpq{p}{q}}{\Gpqf{p}{q}{f}}\vert$.   
\end{corollary}
The theorem and the corollary are illustrated with examples in~\cite{part3}.
Maple worksheets verifying this and other results
from~\cite{part1,part2,part3} can be accessed from~\cite{worksheets}.

\section{Transposition scalar product on spinor spaces}
\label{tpscalarproduct}

In \cite[Ch. 18]{lounesto}, Lounesto discusses scalar products on
$S = \cl_{p,q}f$ for simple Clifford algebras and on $\check{S} =
S \oplus \hat{S} = \cl_{p,q}e$, $e = f+\hat{f}$, for semisimple Clifford
algebras where~$\fg$ denotes the grade involution of $f$. It is well known
that in each case the spinor representation is faithful. Following Lounesto,
we let $\check{\BK}$ be either $\BK$ or $\BK \oplus \hat{\BK}$ and
$\check{S}$ be either $S$ or $S \oplus \hat{S}$ when $\cl_{p,q}$ is simple
or semisimple, respectively. Then, in the simple algebras, the two
$\beta$-scalar products are
\begin{gather}
S \times S \rightarrow \BK, \quad (\psi,\phi) \mapsto 
\begin{cases} \beta_{+}(\psi,\phi) = s_1 \tilde{\psi}\phi \\
              \beta_{-}(\psi,\phi) = s_2 \bar{\psi}\phi 
\end{cases}
\label{eq:betas} 
\end{gather}
whereas in the semisimple algebras they are
\begin{gather}
\chS \times \chS \rightarrow \chK, \quad (\chpsi,\chphi) \mapsto 
\begin{cases} (\beta_{+}(\psi,\phi),\beta_{+}(\psig,\phig)) = (s_1 \tilde{\psi}\phi,s_1 \tilde{\psig}\phig) \\
              (\beta_{-}(\psi,\phi),\beta_{-}(\psig,\phig)) = (s_2 \bar{\psi}\phi,s_2 \bar{\psig}\phig) 
\end{cases}
\label{eq:betass} 
\end{gather}
for $\chpsi=\psi+\psig$ and $\chphi=\phi+\phig$, $\psi,\phi \in S$,
$\psig,\phig\in\hat{S}$, and where $\tilde{\psi},\tilde{\psig}$ (resp.
$\bar{\psi},\bar{\psig}$) denotes reversion (resp. Clifford conjugation) of
$\psi,\psig$. Here $s_1,s_2$ are special monomials in the Clifford algebra
basis $\cb{B}$ which guarantee that the products
$s_1 \tilde{\psi}\phi, \, s_2 \bar{\psi}\phi$, hence also
$s_1 \tilde{\psig}\phig, \, s_2 \bar{\psig}\phig$, belong to
$\BK\cong\hat{\BK}$.\footnote{%
In simple Clifford algebras, the monomials $s_1$ and $s_2$ also satisfy:
(i) $\tilde{f} = s_1 f s_1^{-1}$ and
(ii) $\bar{f} = s_2 f s_2^{-1}$. The identity (i) (resp. (ii)) is also valid
in the semisimple algebras provided $\beta_{+} \not\equiv 0$ (resp.
$\beta_{-} \not\equiv 0$).} 
In fact, the monomials $s_1,s_2$ belong to the chosen transversal of the
stabilizer $\Gpqf{p}{q}{f}$ in $\Gpq{p}{q}$~\cite{part3}. The automorphism
groups of $\beta_{+}$ and $\beta_{-}$ are defined in the simple case as,
respectively, $G_{+}=\{s \in \cl_{p,q} \mid s \tilde{s} =1 \}$ and
$G_{-}=\{s \in \cl_{p,q} \mid s \bar{s} =1 \}$, and as $\fpower{G_{-}}$ and
$\fpower{G_{+}}$ in the semisimple case. They are shown
in~\cite[Tables 1 and 2, p. 236]{lounesto}. 

\subsection{Simple Clifford algebras}
\label{simple}

In Example 3~\cite{part2} it was shown that the transposition scalar product
in $S=\cl_{2,2}f$ is different from each of the two Lounesto's products
whereas Example~4 showed that the transposition product in $S=\cl_{3,0}f$
coincided with $\beta_{+}$. Furthermore, it was remarked that $\tp(\psi)\phi$
always equaled~$\beta_{+}$ for Euclidean signatures $(p,0)$ and~$\beta_{-}$
for anti-Euclidean signatures $(0,q)$. We formalize this in the following
proposition. For all proofs see~\cite{part3}.
\begin{proposition}
Let $\psi,\phi \in S=\cl_{p,q}f$ and $(\psi,\phi) \mapsto \tp(\psi)\phi =
\lambda f$, $\lambda \in \BK,$ be the transposition scalar product. Let
$\beta_{+}$ and $\beta_{-}$ be the scalar products on~$S$ shown
in~(\ref{eq:betas}). Then, there exist monomials $s_1, s_2$ in the transversal
$\ell$ of $\Gpqf{p}{q}{f}$ in $\Gpq{p}{q}$ such that
\begin{gather}
\tp(\psi)\phi
  = \begin{cases} 
      \beta_{+}(\psi,\phi)    = s_1 \tilde{\psi}\phi, \quad \forall \psi,\phi \in \cl_{p,0}f,\\
          \beta_{-}(\psi,\phi)= s_2 \bar{\psi}\phi,   \quad \forall \psi,\phi \in \cl_{0,q}f.
    \end{cases}
    \label{eq:tp=beta}
\end{gather}
\label{prop1} 
\end{proposition}
Let $u \in \cl_{p,q}$ and let $[u]$ be a matrix of $u$ in the spinor
representation $\pi_S$ of $\cl_{p,q}$ realized in the spinor
$(\cl_{p,q},\BK)$-bimodule ${}_{\cl}S_{\BK} \cong \cl_{p,q}f\BK$. Then,
by~\cite[Prop.~5]{part2}, 
\begin{equation}
[\tp(u)] = \begin{cases} 
    [u]^T        & \textit{if $p-q =0,1,2 \bmod 8;$} \\
    [u]^\dagger  & \textit{if $p-q =3,7 \bmod 8;$} \\
    [u]^\ddagger & \textit{if $p-q =4,5,6 \bmod 8;$}
\end{cases}
\label{eq:matdaggers}
\end{equation}
where $T$, $\dagger$, and $\ddagger$ denote, respectively, transposition,
complex Hermitian conjugation, and quaternionic Hermitian conjugation. Thus,
we immediately have:
\begin{proposition}
Let $\Gpqe{p}{q} \subset \cl_{p,q}$ where $\cl_{p,q}$ is a simple Clifford
algebra. Then, $\Gpqe{p}{q}$ is: The orthogonal group $O(N)$ when
$\BK \cong \BR$; the complex unitary group $U(N)$ when $\BK \cong \BC$; or,
the compact symplectic group $Sp(N) = U_\BH(N)$ when
$\BK \cong \BH$.\footnote{%
See Fulton and Harris~\cite{fultonharris} for a definition of the quaternionic
unitary group $U_\BH(N)$. In our notation we follow \textit{loc. cit.}
page 100, `Remark on Notations'.}
That is,
\begin{equation}
\Gpqe{p}{q} = \begin{cases} 
    O(N)  & \textit{if $p-q =0,1,2 \bmod 8;$} \\
    U(N)  & \textit{if $p-q =3,7 \bmod 8;$} \\
    Sp(N) & \textit{if $p-q =4,5,6 \bmod 8;$}
\end{cases}
\label{eq:simpleGpqe}
\end{equation}
where $N=2^k$ and $ k=q-r_{q-p}$. 
\label{propsimpleGpqe}
\end{proposition}

{\small
\begin{table}[t]
\label{tab:t11}
\begin{center}
\renewcommand{\arraystretch}{1.4}
\begin{tabular}{|p{0.4in}|p{0.4in}|p{0.4in}|p{0.41in}|p{0.4in}|p{0.4in}|p{0.4in}|c|}
\multicolumn{8}{c}
{\bf Table 1 (Part 1): Automorphism group $\Gpqe{p}{q}$ of $\tp(\psi)\phi$}\\
\multicolumn{8}{c}
{\bf in simple Clifford algebras $\cl_{p,q} \cong \Mat(2^k,\BR)$}\\
\multicolumn{8}{c}
{$k=q-r_{q-p}$, $p-q \neq 1 \bmod 4, \, p-q = 0,1,2 \bmod 8$}\\\hline
\centering{$(p,q)$} & \centering{$(0,0)$} & \centering{$(1,1)$} & \centering{$(2,0)$} & \centering{$(2,2)$} & \centering{$(3,1)$} & \centering{$(3,3)$} & $(0,6)$ \\\hline
\centering{$\Gpqe{p}{q}$} & \centering{$O(1)$} & \centering{$O(2)$} & \centering{$\boxed{O(2)}$} &\centering{$O(4)$} & \centering{$O(4)$} & \centering{$O(8)$} & 
$\boxed{\boxed{O(8)}}$ \\\hline
\end{tabular}
\end{center}
\end{table}
}

{\small
\begin{table}[t]
\label{tab:t12}
\begin{center}
\renewcommand{\arraystretch}{1.4}
\begin{tabular}{|p{0.42in}|p{0.42in}|p{0.47in}|p{0.47in}|p{0.62in}|p{0.47in}|c|}
\multicolumn{7}{c}
{\bf Table 1 (Part 2): Automorphism group $\Gpqe{p}{q}$ of $\tp(\psi)\phi$}\\
\multicolumn{7}{c}
{\bf in simple Clifford algebras $\cl_{p,q} \cong \Mat(2^k,\BR)$}\\
\multicolumn{7}{c}
{$k=q-r_{q-p}$, $p-q \neq 1 \bmod 4, \, p-q = 0,1,2 \bmod 8$}\\\hline
\centering{$(p,q)$} & \centering{$(4,2)$} & \centering{$(5,3)$} & \centering{$(1,7)$} & \centering{$(0,8)$} & \centering{$(4,4)$} & $(8,0)$\\\hline
\centering{$\Gpqe{p}{q}$} & \centering{$O(8)$} & \centering{$O(16)$} & \centering{$O(16)$} & \centering{$\boxed{\boxed{O(16)}}$} & \centering{$O(16)$} & $\boxed{O(16)}$ \\\hline
\end{tabular}
\end{center}
\end{table}
}

{\small
\begin{table}[t]
\label{tab:t21}
\begin{center}
\renewcommand{\arraystretch}{1.4}
\begin{tabular}{|c|c|c|c|c|c|c|c|c|}
\multicolumn{9}{c}
{\bf Table 2 (Part 1): Automorphism group $\Gpqe{p}{q}$ of $\tp(\psi)\phi$}\\
\multicolumn{9}{c}
{\bf in simple Clifford algebras $\cl_{p,q} \cong \Mat(2^k,\BC)$}\\
\multicolumn{9}{c}
{$k=q-r_{q-p}$, $p-q \neq 1 \bmod 4, \, p-q = 3,7 \bmod 8$}\\\hline
$(p,q)$ & $(0,1)$ & $(1,2)$ & $(3,0)$ & $(2,3)$ & $(0,5)$ & $(4,1)$ & $(1,6)$ & $(7,0)$ \\\hline
$\Gpqe{p}{q}$ & $U(1)$ & $U(2)$ & $\boxed{U(2)}$ & $U(4)$ & $\boxed{\boxed{U(4)}}$ & $U(4)$ & $U(8)$ & $\boxed{U(8)}$ \\\hline
\end{tabular}
\end{center}
\end{table}
}

{\small
\begin{table}[t]
\label{tab:t22}
\begin{center}
\renewcommand{\arraystretch}{1.4}
\begin{tabular}{|c|c|c|c|c|c|c|c|}
\multicolumn{8}{c}
{\bf Table 2 (Part 2): Automorphism group $\Gpqe{p}{q}$ of $\tp(\psi)\phi$}\\
\multicolumn{8}{c}
{\bf in simple Clifford algebras $\cl_{p,q} \cong \Mat(2^k,\BC)$}\\
\multicolumn{8}{c}
{$k=q-r_{q-p}$, $p-q \neq 1 \bmod 4, \, p-q = 3,7 \bmod 8$}\\\hline
$(p,q)$ & $(5,2)$ & $(3,4)$ & $(4,5)$ & $(6,3)$ & $(2,7)$ & $(0,9)$ & $(8,1)$\\\hline
$\Gpqe{p}{q}$ & $U(8)$ & $U(8)$ & $U(16)$ & $U(16)$ & $U(16)$ & $\boxed{\boxed{U(16)}}$ & $U(16)$ \\\hline
\end{tabular}
\end{center}
\end{table}
}

{\small
\begin{table}[t]
\label{tab:t31}
\begin{center}
\renewcommand{\arraystretch}{1.4}
\begin{tabular}{|p{0.41in}|p{0.6in}|p{0.6in}|p{0.5in}|p{0.43in}|p{0.43in}|c|}
\multicolumn{7}{c}
{\bf Table 3 (Part 1): Automorphism group $\Gpqe{p}{q}$ of $\tp(\psi)\phi$}\\
\multicolumn{7}{c}
{\bf in simple Clifford algebras $\cl_{p,q} \cong \Mat(2^k,\BH)$}\\
\multicolumn{7}{c}
{$k=q-r_{q-p}$, $p-q \neq 1 \bmod 4, \, p-q = 4,5,6 \bmod 8$}\\\hline
\centering{$(p,q)$} & \centering{$(0,2)$} & \centering{$(0,4)$} & \centering{$(4,0)$} & \centering{$(1,3)$} & \centering{$(2,4)$} & $(6,0)$\\\hline
\centering{$\Gpqe{p}{q}$} & \centering{$\boxed{\boxed{Sp(1)}}$} & 
\centering{$\boxed{\boxed{Sp(2)}}$} & \centering{$\boxed{Sp(2)}$} & \centering{$Sp(2)$} & \centering{$Sp(4)$} & $\boxed{Sp(4)}$\\\hline
\end{tabular}
\end{center}
\end{table}
}%

{\small
\begin{table}[t]
\label{tab:t32}
\begin{center}
\renewcommand{\arraystretch}{1.4}
\begin{tabular}{|p{0.42in}|p{0.48in}|p{0.48in}|p{0.48in}|p{0.48in}|p{0.48in}|c|}
\multicolumn{7}{c}
{\bf Table 3 (Part 2): Automorphism group $\Gpqe{p}{q}$ of $\tp(\psi)\phi$}\\
\multicolumn{7}{c}
{\bf   in simple Clifford algebras $\cl_{p,q} \cong \Mat(2^k,\BH)$}\\
\multicolumn{7}{c}
{$k=q-r_{q-p},$ $p-q \neq 1 \bmod 4, \, p-q = 4,5,6 \bmod 8$}\\\hline
\centering{$(p,q)$} & \centering{$(1,5)$} & \centering{$(5,1)$} & \centering{$(6,2)$} & \centering{$(7,1)$} & \centering{$(2,6)$} & $(3,5)$\\\hline
\centering{$\Gpqe{p}{q}$} & \centering{$Sp(4)$} & \centering{$Sp(4)$} & \centering{$Sp(8)$} & \centering{$Sp(8)$} & \centering{$Sp(8)$} & $Sp(8)$\\\hline
\end{tabular}
\end{center}
\end{table}
}

The scalar product $\tp(\psi)\phi$ was computed with
$\mathtt{CLIFFORD}$~\cite{clifford3} for all signatures $(p,q)$,
$p+q \leq 9$~\cite{worksheets}. Observe that as expected, in Euclidean (resp.
anti-Euclidean) signatures $(p,0)$ (resp. $(0,q)$) the group $\Gpqe{p}{0}$
(resp. $\Gpqe{0}{q}$) coincides with the corresponding automorphism group of
the scalar product~$\beta_{+}$ (resp. $\beta_{-}$) listed 
in~\cite[Table 1, p. 236]{lounesto} (resp.~\cite[Table 2, p. 236]{lounesto}).
This is indicated by a single (resp. double) box around the group symbol in
Tables~1--5. For example, in Table~1, for the Euclidean signature $(2,0)$, we
show $\Gpqe{2}{0}$ as $\boxed{O(2)}$ like for $\beta_{-}$ whereas for the
anti-Euclidean signature $(0,6)$, we show $\Gpqe{0}{6}$ as
$\boxed{\boxed{O(8)}}$ like for $\beta_{+}$.

For simple Clifford algebras, the automorphism groups $\Gpqe{p}{q}$ are
displayed in Tables~1,~2, and~3. In each case the form is positive definite
and non-degenerate. Also, unlike in the case of the forms $\beta_{+}$ and
$\beta_{-}$, there is no need for the extra monomial factor like $s_1,s_2$
in~(\ref{eq:betas}) (and~(\ref{eq:betass})) to guarantee that the product
$\tp(\psi)\phi$ belongs to $\BK$ since this is always the
case~\cite{part1, part2}. Recall that the only role of the monomials $s_1$
and $s_2$ is to permute entries of the spinors $\tilde{\psi}\phi$
and $\bar{\psi}\phi$ to assure that $\beta_{+}(\psi,\phi)$ and
$\beta_{-}(\psi,\phi)$ belong to the (skew) field~$\BK$. That is, more
precisely, that $\beta_{+}(\psi,\phi)$ and $\beta_{-}(\psi,\phi)$ have the
form $\lambda f = f \lambda$ for some $\lambda$ in $\BK$. The idempotent $f$
in the spinor basis in $S$ corresponds uniquely to the identity coset
$\Gpqf{p}{q}{f}$ in the quotient group $\Gpq{p}{q}/\Gpqf{p}{q}{f}$. Based
on~\cite[Prop. 2]{part2} we know that since the vee group $\Gpq{p}{q}$
permutes entries of any spinor~$\psi$, the monomials $s_1$ and $s_2$ belong
to the transversal of the stabilizer
$\Gpqf{p}{q}{f} \lhd \Gpq{p}{q}$~\cite[Cor. 2]{part2}.\footnote{%
In \cite[Page 233]{lounesto}, Lounesto states correctly that ``the element
$s$ can be chosen from the standard basis of $\cl_{p,q}$." In fact, one can
restrict the search for $s$ to the transversal of the stabilizer
$\Gpqf{p}{q}{f}$ in $\Gpq{p}{q}$ which has a much smaller size
$2^{q-r_{q-p}}$ compared to the size $2^{p+q}$ of the Clifford basis.}

One more difference between the scalar products $\beta_{+}$ and $\beta_{-}$,
and the transposition product $\tp(\psi)\phi$ is that in some signatures one
of the former products may be identically zero whereas the transposition
product is never identically zero. The signatures $(p,q)$ in which one of the
products $\beta_{+}$ or $\beta_{-}$ is identically zero can be easily found
in~\cite[Tables 1 and 2, p. 236]{lounesto} as the automorphism group of the
product is then a general linear group. 

\subsection{Semisimple Clifford algebras}
\label{semisimple}
Faithful spinor representation of a semisimple Clifford algebra $\cl_{p,q}$
($p-q =1 \bmod 4$) is realized in a left ideal $\check{S} = S\oplus \hat{S}
= \cl_{p,q}e$ where $e = f + \hat{f}$ for any  primitive idempotent~$f$.
Recall that $\hat{\phantom{u}}$ denotes the grade involution of
$u \in \cl_{p,q}$. We refer to~\cite[pp. 232--236]{lounesto} for some of the
concepts. In particular, $S = \cl_{p,q}f$ and $\Sg = \cl_{p,q}\fg$. Thus,
every spinor $\check{\psi} \in \check{S}$ has unique components $\psi \in S$
and $\psig \in \Sg$. We refer to the elements $\check{\psi} \in \check{S}$
as ``spinors" whereas to its components $\psi \in S$ and $\psig \in \Sg$ we
refer as ``$\frac12$-spinors". 

For the semisimple Clifford algebras $\cl_{p,q}$, we will view spinors
$\check{\psi} \in \check{S}= S \oplus \hat{S}$ as ordered pairs
$(\psi,\psi_g) \in S \times \hat{S}$ when 
$\check{\psi} = \psi + \psi_g$. Likewise, we will view elements
$\check{\lambda}$ in the double fields $\check{\BK} = \BK \oplus \hat{\BK}$ as
ordered pairs $(\lambda,\lambda_g)\in \BK \times \hat{\BK}$ when
$\check{\lambda} = \lambda + \lambda_g$. As before, $\BK = f \cl_{p,q}f$ while
$\hat{\BK} = \hat{f} \cl_{p,q}\hat{f}$. Recall that 
$\check{\BK} \cong \fpower{\BR} \stackrel{\mbox{\scriptsize{def}}}{=} 
\BR \oplus \BR$ or
$\check{\BK} \cong \fpower{\BH} \stackrel{\mbox{\scriptsize{def}}}{=} 
\BH \oplus \BH$ when, respectively, $p-q = 1 \bmod 8$, or $p-q = 5 \bmod 8$.

In this section we classify automorphism groups of the transposition scalar
product
\begin{gather}
  \check{S} \times \check{S} \rightarrow \check{\BK}, \quad (\check{\psi},\check{\phi}) \mapsto \tp(\check{\psi},\check{\phi}) \stackrel{\mbox{\scriptsize{def}}}{=} 
  (\tp(\psi)\phi, \tp(\psig)\phig) \in \check{\BK}
\label{eq:tsprod}
\end{gather}
when $\chpsi = \psi+\psig$ and $\chphi = \phi+\phig$.

\begin{proposition}
Let $\Gpqe{p}{q} \subset \cl_{p,q}$ where $\cl_{p,q}$ is a semisimple Clifford
algebra. Then, $\Gpqe{p}{q}$ is: The double orthogonal group 
$\fpower{O(N)} \stackrel{\mbox{\scriptsize{def}}}{=} O(N) \times O(N)$ when
$\check{\BK} \cong \fpower{\BR}$ or the double compact symplectic group 
$\fpower{Sp(N)} \stackrel{\mbox{\scriptsize{def}}}{=} Sp(N) \times Sp(N)$ when
$\check{\BK} \cong \fpower{\BH}$.\footnote{%
Recall that $Sp(N)=U_{\BH}(N)$ where $U_{\BH}(N)$ is the quaternionic
unitary group~\cite{fultonharris}.}
That is,
\begin{equation}
  \Gpqe{p}{q}
    = \begin{cases} 
        \fpower{O(N)}  = O(N) \times O(N)   & \textit{when $p-q =1 \bmod 8;$} \\
        \fpower{Sp(N)} = Sp(N) \times Sp(N) & \textit{when $p-q =5 \bmod 8;$}
      \end{cases}
\label{eq:semisimpleGpqe}
\end{equation}
where $N=2^{k-1}$ and $k=q-r_{q-p}$. 
\label{propsemisimpleGpqe}
\end{proposition}

The automorphism groups $\Gpqe{p}{q}$ for semisimple Clifford algebras
$\cl_{p,q}$ for $p+q\leq 9$ are shown in Tables~4 and~5. All results in these
tables, like in Tables 1, 2, and 3, have been verified with
$\mathtt{CLIFFORD}$~\cite{clifford3} and the corresponding Maple worksheets
are posted at~\cite{worksheets}. 

{\small
\begin{table}[t]
\label{tab:t4}
\begin{center}
\renewcommand{\arraystretch}{1.4}
\begin{tabular}{|p{0.3in}|p{0.46in}|p{0.33in}|p{0.33in}|p{0.56in}|p{0.33in}|p{0.43in}|p{0.43in}|c|}
\multicolumn{9}{c}
{\bf Table 4: Automorphism group $\Gpqe{p}{q}$ of $\tp(\psi)\phi$}\\
\multicolumn{9}{c}
{\bf in semisimple Clifford algebras $\cl_{p,q} \cong \fpower{\Mat(2^{k-1},\BR)}$}\\
\multicolumn{9}{c}
{$k=q-r_{q-p}, $ $p-q = 1 \bmod 4, \, p-q = 1 \bmod 8$}\\\hline
\centering{$(p,q)$} & \centering{$(1,0)$} & \centering{$(2,1)$} & \centering{$(3,2)$} & \centering{$(0,7)$} & \centering{$(4,3)$} & \centering{$(1,8)$} & \centering{$(5,4)$} & $(9,0)$\\\hline
\centering{$\Gpqe{p}{q}$} & \centering{$\boxed{\fpower{O(1)}}$} & \centering{$\fpower{O(2)}$} & \centering{$\fpower{O(4)}$} & \centering{$\boxed{\boxed{\fpower{O(8)}}}$} & 
\centering{$\fpower{O(8)}$} & \centering{$\fpower{O(16)}$} & \centering{$\fpower{O(16)}$} & $\boxed{\fpower{O(16)}}$ \\\hline 
\end{tabular}
\end{center}
\end{table}
}%

{\small
\begin{table}[t]
\label{tab:t5}
\begin{center}
\renewcommand{\arraystretch}{1.4}
\begin{tabular}{|c|c|c|c|c|c|c|c|}
\multicolumn{8}{c}
{\bf Table 5: Automorphism group $\Gpqe{p}{q}$ of $\tp(\psi)\phi$}\\
\multicolumn{8}{c}
{\bf in semisimple Clifford algebras $\cl_{p,q} \cong \fpower{\Mat(2^{k-1},\BH)}$}\\
\multicolumn{8}{c}
{$k=q-r_{q-p},$ $p-q = 5 \bmod 4, \, p-q =  \bmod 8$}\\\hline
$(p,q)$ & $(0,3)$ & $(1,4)$ & $(5,0)$ & $(2,5)$ & $(6,1)$ & $(3,6)$ & $(7,2)$\\\hline
$\Gpqe{p}{q}$ & $\boxed{\boxed{\fpower{Sp(1)}}}$ & $\fpower{Sp(2)}$ & 
$\boxed{\fpower{Sp(2)}}$ & $\fpower{Sp(4)}$ & $\fpower{Sp(4)}$ & $\fpower{Sp(8)}$ & 
$\fpower{Sp(8)}$ \\\hline
\end{tabular}
\end{center}
\end{table}
}%

\section{Conclusions}
\label{conc}

The transposition map $\tp$ allowed us to define a new transposition scalar
product on spinor spaces. Only in the Euclidean and anti-Euclidean signatures,
this scalar product is identical to the two known spinor scalar products
$\beta_{+}$ and $\beta_{-}$ which use, respectively, the reversion and the
conjugation and it is different in all other signatures. This new product is
never identically zero and it does not require extra monomial factor to assure
it is $\BK$- or $\check{\BK}$-valued. This is because the $\tp$ maps any
spinor space to its dual. Then, we have identified the automorphism groups
$\Gpqe{p}{q}$ of this new product in Tables 1--5 for $p+q=n\leq 9$. The
classification is complete and sufficient due to the $\bmod\, 8$ periodicity.

We have observed the important role played by the idempotent group
$\Tpqf{p}{q}{f}$ and the field group $\Kpqf{p}{q}{f}$ as normal subgroups in
the stabilizer group $\Gpqf{p}{q}{f}$ of the primitive idempotent $f$ and
their coset spaces $\Gpq{p}{q}/\Tpqf{p}{q}{f}$,
$\Gpqf{p}{q}{f}/\Tpqf{p}{q}{f}$, and $\Gpq{p}{q}/\Gpqf{p}{q}{f}$ in relation
to the spinor representation of $\cl_{p,q}$. These subgroups allow to
construct very effectively non-canonical transversals and hence basis elements
of the spinor spaces and the (skew double) field underlying the spinor space.
This approach to the spinor representation of $\cl_{p,q}$ based on the
stabilizer $\Gpqf{p}{q}{f}$ of $f$ leads to a realization that the Clifford
algebras can be viewed as a twisted group ring $\BR^t[(\BZ_2)^n]$. In
particular, we have observed that our transposition $\tp$ is then a `star map'
of $\BR^t[(\BZ_2)^n]$~\cite{passman} which on a general twisted group ring
$*:K^t[G] \rightarrow K^t[G]$ is defined as 
$$
  \left(\sum a_x \bar{x}\right)^* = \sum a_x \bar{x}^{-1}. 
$$
This is because we recall properties of the transposition anti-involution
$\tp$, and, in particular, its action $\tp(m)=m^{-1}$ on a monomial $m$ in the
Grassmann basis $\cb{B}$ which is, as we see now, identical to the action
$*(m)=m^{-1}$ on every $m \in \cb{B}$. For a Hopf algebraic discussion of
Clifford algebras as twisted group algebras,
see~\cite{albuquerquemajid,morier-genoudovsienko} and references therein.



\begin{thebibliography}{00}

\bibitem{ablamowicz1996}
R.~Ab\l amowicz, \textit{Clifford Algebra Computations with Maple}, in ``Clifford (Geometric) Algebras with Applications in Physics, Mathematics, and Engineering'', W.~E.~Baylis, (Ed.) (Birkh\"{a}user, Boston, 1996) 463--502
\bibitem{ablamowicz1998}
R.~Ab\l amowicz, \textit{Spinor Representations of Clifford algebras: A Symbolic Approach}, Computer Physics Communications Thematic Issue - Computer Algebra in Physics Research \textbf{115}, No. 2--3 (1998) 510--535
\bibitem{ablamowicz2005}  R.~Ab\l amowicz, \textit{Computations with Clifford and Grassmann Algebras}. Adv. in Appl. Clifford Algebras, \textbf{19}, No. 3--4 (2009) 499--545
\bibitem{ablamowicz2009a} R.~Ab\l amowicz, \textit{Computation of Non-Commutative Gr\"{o}bner Bases in Grassmann and Clifford Algebras}. Adv. Applied Clifford Algebras 20 No. 3--4 (2010), 447--476
\bibitem{GfG} R.~Ab\l amowicz and B.~Fauser, \textit{$\mathtt{GfG}$ - Groebner for Grassmann - A Maple 12 Package for Groebner Bases in Grassmann Algebras}, 
{\texttt \small http://math.tntech.edu/rafal/GfG12/} (2010)
\bibitem{part1}
R.~Ab\l amowicz and B.~Fauser, \textit{On the Transposition Anti-Involution in Real Clifford Algebras I: The Transposition Map}. Linear and 
Multilinear Algebra, Volume 59, Issue 12, December 2011, 1331--1358
\bibitem{part2}
R.~Ab\l amowicz and B.~Fauser, \textit{On the Transposition Anti-Involution in Real Clifford Algebras II: Stabilizer Groups of Primitive Idempotents}. Linear and Multilinear Algebra, Volume 59, Issue 12, December 2011, 1359--1381
\bibitem{part3} R.~Ab\l amowicz and B.~Fauser, \textit{On the Transposition Anti-Involution in Real Clifford Algebras III: The Automorphism Group of the Transposition Scalar Product on Spinor Spaces} (to appear in Linear and Multilinear Algebra, DOI: 10.1080/03081087.2011.624093)
\bibitem{clifford3}
R.~Ab\l amowicz and B.~Fauser, $\mathtt{CLIFFORD}$ with $\mathtt{Bigebra}$ -- A Maple Package for Computations with Clifford and Grassmann Algebras, 
{\texttt \small http://math.tntech.edu/rafal/} (\copyright 1996--2011)
\bibitem{worksheets}
R.~Ab\l amowicz and B.~Fauser, Maple worksheets created with $\mathtt{CLIFFORD}$ for a verification of results presented in this paper and in~\cite{part1,part2}, 
{\texttt \small http://math.tntech.edu/rafal/publications.html} (\copyright 2011)
\bibitem{albuquerquemajid}
H.~Albuquerque and S.~Majid, \textit{Clifford Algebras Obtained by Twisting of Group Algebras},
J. of Pure and Appl. Algebra \textbf{171}, (2002) 133--148
\bibitem{fultonharris}
W.~Fulton and J.~Harris, \textit{Representation Theory: A First Course}, (Springer, New York, 1991)
\bibitem{hahn}
A.~J.~Hahn, \textit{Quadratic Algebras, Clifford Algebras, and Arithmetic Witt Groups},
(Undergraduate Texts in Mathematics) (Springer-Verlag, New York, 1994) 
\bibitem{lounesto}
P.~Lounesto,
\textit{Clifford Algebras and Spinors}, 2nd ed. (Cambridge University Press, Cambridge, 2001)
\bibitem{morier-genoudovsienko}
S.~Morier-Genoud and V. Ovsienko, \textit{Simple graded commutative algebras}, J. of Algebra \textbf{323} (2010) 1649--1664
\bibitem{oziewicz1997} Z.~Oziewicz, 
\textit{Clifford hopf-gebra and bi-universal hopf-gebra},
\textit{Czechoslovak Journal of Physics} Volume 47, Number 12, 1267-1274, DOI: 10.1023/A:1022833801475
\bibitem{passman}
D.~S.~Passman, 
\textit{The Algebraic Structure of Group Rings}, (Robert E. Krieger Publishing Company, Malabar, Florida, 1985)
\bibitem{porteous}
I.~R.~Porteous, \textit{Clifford Algebras and the Classical Groups}, (Cambridge University Press, Cambridge, 1995)
\bibitem{rotman}
J.J.~Rotman,
\textit{Advanced Modern Algebra}, Revised Printing, (Prentice Hall, Upper Saddle River, 2002)
\bibitem{salingaros1}
N.~Salingaros, 
\textit{Realization, extension, and classification of certain physically important groups and algebras}, J. Math. Phys. \textbf{22} (1981) 226--232
\bibitem{salingaros2}
N.~Salingaros, \textit{On the classification of Clifford algebras and their relation to spinors in $n$ dimensions}, J. Math. Phys. \textbf{23} (1) (1982) 1--7
\bibitem{salingaros3}
N.~Salingaros, \textit{The relationship between finite groups and Clifford algebras}, J. Math. Phys. \textbf{25} (1984) 738--742
%
\end{thebibliography}
\end{document}